\documentclass[a4paper,12pt]{article}

\title{
Simple Proofs for the Derivative Estimates of the Holomorphic Motion near Two Boundary Points of the Mandelbrot Set\\
}
\author{
Yi-Chiuan Chen 
and
Tomoki Kawahira
}

\usepackage{enumerate}
\usepackage{graphicx}
\usepackage{amsmath}
\usepackage{amssymb}
\usepackage{amsthm}

\theoremstyle{plain}
\newtheorem{thm}{Theorem}[section]
\newtheorem{prop}[thm]{Proposition}
\newtheorem{lem}[thm]{Lemma}
\newtheorem{cor}[thm]{Corollary}
\newtheorem{remark}[thm]{Remark}
\theoremstyle{remark}
\theoremstyle{definition}

\newcommand{\C}{\mathbb{C}}

\newcommand{\M}{\mathbb{M}}
\newcommand{\T}{\mathbb{T}}

\newcommand{\chat}{{\hat{c}}}
\newcommand{\abs}[1]{{\left| #1 \right|}}
\newcommand{\paren}[1]{{\left( #1 \right)}}
\newcommand{\brac}[1]{{\left\{ #1 \right\}}}

\newcommand{\cR}{{\mathcal{R}}}

\newcommand{\QED}{\hfill $\blacksquare$}

\begin{document}

\maketitle

\begin{abstract}
For the complex quadratic family $q_c:z\mapsto z^2+c$,
it is known that every point in the Julia set $J(q_c)$ 
moves holomorphically on $c$ except 
at the boundary points of the Mandelbrot set. 
In this note, we present short proofs of the following derivative estimates of the motions near the boundary points $1/4$ and $-2$:
for each $z = z(c)$ in the Julia set, the derivative $dz(c)/dc$ 
is uniformly $O(1/\sqrt{1/4-c})$ when real $c\nearrow 1/4$;
and is uniformly $O(1/\sqrt{-2-c})$ when real $c\nearrow -2$.
These estimates of the derivative imply Hausdorff convergence of the Julia set $J(q_c)$ when $c$ approaches these boundary points.
In particular, the Hausdorff distance between $J(q_c)$ with $0\le c<1/4$
and $J(q_{1/4})$ is exactly $\sqrt{1/4-c}$.
\end{abstract}


\section{Introduction}

For the family of quadratic maps $q_{c}$ of $\C$, $z\mapsto z^2+c$, with $c$ a complex number not locating on the boundary of the Mandelbrot set $\M$, 
it is well-known that every point in the Julia set $J(q_c)$ moves holomorphically with respect to $c$, i.e. the {\it holomorphic motion} \cite{L, MSS}.  
Note that $J(q_c)$ is a Cantor set when $c\not\in\M$ and 
is connected when $c\in\M$.
A parameter $c$ is called {\it hyperbolic} if the orbit of the origin 
accumulates on an attracting cycle or diverges to infinity. 
A {\it hyperbolic component} is a connected component of $\C-\partial \M$ 
containing hyperbolic parameters.
It is conjectured that the set $\C-\partial \M$ 
consists of only the hyperbolic parameters. 
See \cite[Expos\'e I]{DH} for example.

Let $\mathbb{D}$ denote the open disk of radius one centered at the origin.
 There is a biholomorphic function $\Phi$ from $\overline{\mathbb{C}}- \M$ to $\overline{\mathbb{C}}- \overline{\mathbb{D}}$ with which the set
 \[ \mathcal{R}(\theta):=\{\Phi^{-1}(r e^{i2\pi\theta})|~ 1< r\le\infty\}
\]
 is defined and called the {\it parameter ray} of angle $\theta\in\T=\mathbb{R}/\mathbb{Z}$ of the Mandelbrot set $\M$. 
Given $\theta$, if $\lim_{r\searrow 1}\Phi^{-1}(re^{i2\pi\theta})$ exists, then this limit is called the {\it landing point} of the parameter ray $\mathcal{R}(\theta)$.
A parameter $\hat{c}$ in $\partial \M$ is called {\it semi-hyperbolic} 
if the critical point is non-recurrent and belongs to the Julia set \cite{CJY}. A typical example is a {\it Misiurewicz} parameter,
that is, for which the critical point eventually lands on a repelling periodic point. 
The set consisting of semi-hyperbolic parameters is dense with Hausdorff dimension 2 in $\partial\M$ \cite{Shi}.
For each semi-hyperbolic parameter $\chat \in \partial \M$, 
there exists at least one parameter ray $\cR(\theta)$ landing at $\chat$. 
(See \cite[Theorem 2]{D2}.) 

In a recent paper \cite{CK}, we proved the following result concerning the estimate for the derivative of the holomorphic motion.
\begin{thm}  \label{thm_Main_estimate}
Let $\hat{c} \in \partial \M$ be a semi-hyperbolic parameter 
that is the landing point of $\mathcal{R}(\theta)$.
Then there exists a constant $K>0$ that depends only on $\hat{c}$
such that for any $c\in\mathcal{R}(\theta)$ 
sufficiently close to $\chat$ and any $z = z(c)\in J(q_c)$, 
the point $z(c)$ moves holomorphically with 
\begin{equation}
\abs{\frac{dz(c)}{dc}}
\le 
\frac{K}{\sqrt{\abs{ c-\hat{c}}}}. \label{Main_estimate}
\end{equation}
\end{thm}
This result enables us to obtain one-sided H\"older continuity of the holomorphic motion
along the parameter ray (i.e. the holomorphic motion lands). More precisely, 
let $\hat{c}\in \partial \M$ be a semi-hyperbolic parameter 
that is a landing point of $\mathcal{R}(\theta)$, 
and let 
$c=c(r):=\Phi^{-1}(re^{i2\pi\theta})$ with  $r \in (1,2]$.
Then for any $z(c(2))$ in $J(q_{c(2)})$, 
the improper integral 
$$
z(\hat{c}) := z(c(2))+ \lim_{\delta \searrow 0} \int_{2}^{1+\delta} \frac{dz(c)}{dc} ~\frac{dc(r)}{dr}~dr
$$
exists in the Julia set $J(q_{\chat})$. 
In particular, $z(c)$ is uniformly 
one-sided H\"older continuous of exponent $1/2$ at $c = \hat{c}$ along $\mathcal{R}(\theta)$: There exists a constant $K'$ depending only on $\chat$ such that 
\begin{equation}
|z(c)-z(\chat)| \le K' \sqrt{|c-\chat|} \label{eq2}
\end{equation}
for any $c =c(r) \in \mathcal{R}(\theta)$ with $1<r\le 2$.

The primary aims of this paper are two-folds. The first is to show that the same estimate as \eqref{Main_estimate} holds when $c$ approaches $1/4$, which is a parabolic parameter, along the real line in the interior of the Mandelbrot set, namely
\begin{thm} \label{Main_thm_2}
For  any point $z = z(c)$ in the Julia set $J(q_c)$, we have that $z(c)$ moves holomorphically with derivative
\begin{equation}
\abs{\frac{dz(c)}{dc}}
 \le \frac{1}{2\sqrt{1/4-c}} \label{estimate2} 
\end{equation}
as $c \nearrow 1/4$ along the real line in the interior of $\M$.
\end{thm} 

\begin{figure}[htbp]
\begin{center}
\includegraphics[width=.95\textwidth]{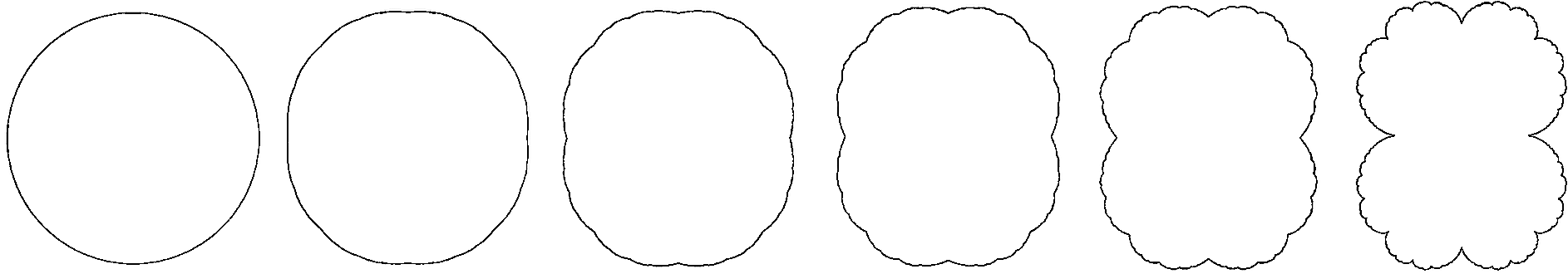}
\includegraphics[width=.55\textwidth]{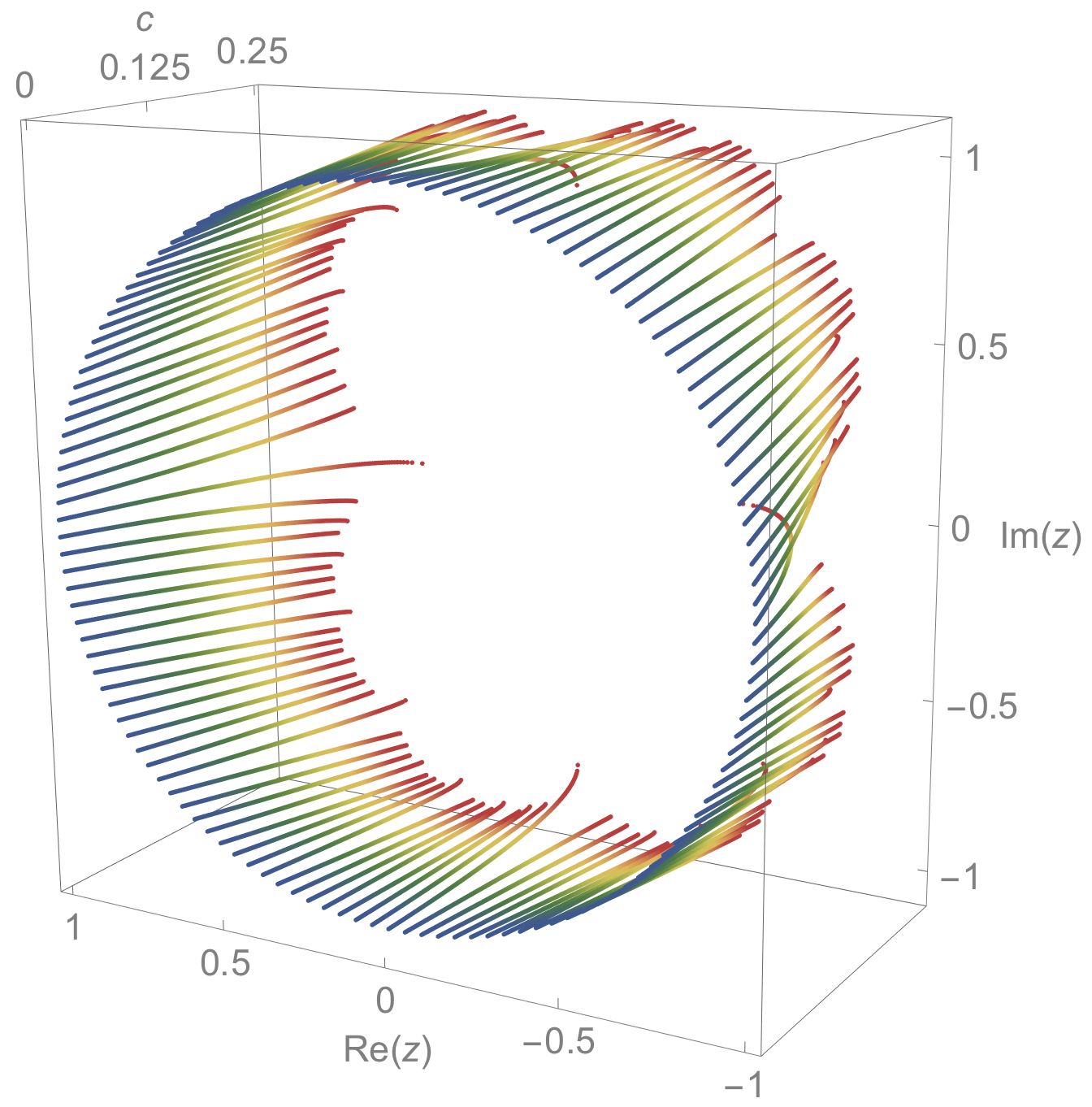}
\end{center}
\caption{Top: The Julia set $J(q_c)$ for $c= k/20~ (k=0, 1, \cdots, 5)$. 
Bottom: Real analytic motion of the preimages of the repelling fixed point for $0 \le c <1/4$.}
\label{fig_cauliflowers}
\end{figure}

The second is to present a simple proof of Theorem \ref{thm_Main_estimate} for the case $\chat=-2$. (See Remark  \ref{rk:thm13}.) As a matter of fact, what we present is a proof of Theorem \ref{thm_Main_estimate} for the logistic map $f_\mu:\C\to\C$, $z\mapsto \mu z(1-z)$, with the semi-hyperbolic (Misiurewicz) parameter  $\mu=4$ case, stated as follows.

\begin{thm}\label{Main_thm}
For  $z = z(\mu)$ in the Julia set $J(f_\mu)$, the point
$z(\mu)$ moves holomorphically with derivative
$$
\abs{\frac{dz(\mu)}{d\mu}}
 =  O\paren{\frac{1}{\sqrt{\mu-4}}} 
$$
as $\mu \searrow 4$ along the real line.
\end{thm}

\begin{figure}[htbp]
\begin{center}
\includegraphics[width=.55\textwidth]{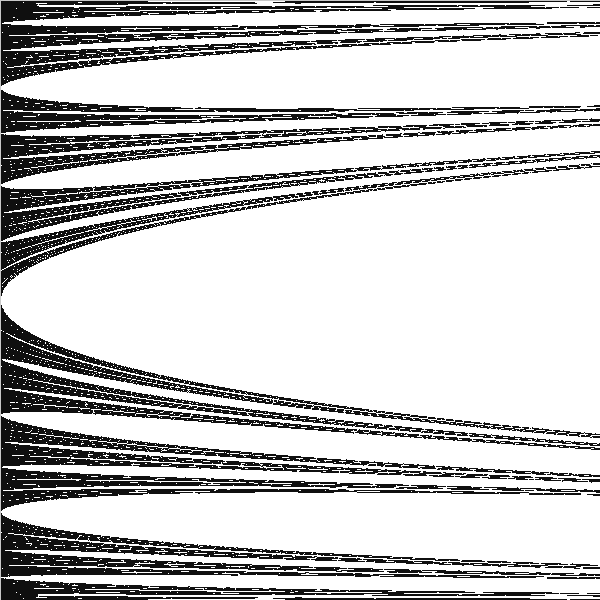}
\end{center}
\caption{Real analytic motion of the Julia set $J(f_\mu)$ for 
$\mu \searrow 4$. See Figure 2 of \cite{CK} for the corresponding motion for $q_c$ for $c \nearrow -2$.}
\label{fig_mu24}
\end{figure}

The proof of Theorem \ref{thm_Main_estimate} in \cite{CK} relies on the hyperbolic metric, the hyperbolicity of the $\omega$-limit set of a semi-hyperbolic parameter $\chat$, the John condition on $\mathbb{C}-J(q_\chat)$, and the asymptotic similarity between $J(q_{\hat{c}})$ and $\M$ at $\hat{c}$. To the best of our knowledge, the difficulty of the proof for parameter $\chat=-2$ is essentially the same as the general semi-hyperbolic parameter case. Moreover, the difficulty of proof remains unchanged even if we consider quadratic maps in the form of logistic maps for general semi-hyperbolic parameter case. To our surprise,   however, we find that for the logistic map with $\mu=4$ case, the difficulty of the proof can be substantially reduced.  We only need a singular metric, and the proof is very straightforward. This is one of the motivations of this paper.

\begin{remark} \rm \label{rk:thm13}
 The logistic map $f_\mu$ is affinely conjugate to $q_c$ via the conjugacy  
\begin{equation}
G(\cdot, \mu):z \mapsto w=-\mu z+\frac{\mu}{2} \qquad (\mbox{so that}~ G(\cdot, \mu)\circ f_\mu\circ G(\cdot, \mu)^{-1}=q_c) \label{conjugacy}  
\end{equation}
with 
\begin{equation}
c=\frac{\mu (2-\mu)}{4} \qquad (\mu\not= 0). \label{cmu}  
\end{equation}
Fix a point $z_0\in J(f_\mu)$, and let $w_0=G(z_0, \mu)\in J(q_c)$. When $\mu$ varies, let $z(\mu)$ be the holomorphic motion for $f_\mu$, and $w(c)$ be the corresponding holomorphic motion for $q_c$ via the above relation \eqref{cmu}. From the conjugacy \eqref{conjugacy}, the derivative $dw(c)/dc$ can be obtained from $dz(\mu)/d\mu$:
\begin{eqnarray}
  \frac{dw(c)}{dc} &=& \frac{\partial G}{\partial z}\cdot\frac{dz}{d\mu}\cdot\frac{d\mu}{dc}+\frac{\partial G}{\partial\mu}\cdot\frac{d\mu}{dc} \nonumber \\
                          &=& \left(-\mu\cdot\frac{dz}{d\mu}-z+\frac{1}{2}\right)\cdot\frac{2}{1-\mu} \qquad (\mu\not=0~\mbox{or}~1).  \label{dwdc}
\end{eqnarray}
Also, we have
\begin{equation}
  \frac{dz(\mu)}{d\mu} = \frac{\mu-1}{2\mu}\cdot\frac{dw}{dc}+\frac{w}{\mu^2}  \qquad (\mu\not=0~\mbox{or}~1).  \label{dzdmu}
\end{equation}
 Theorem \ref{thm_Main_estimate} for  $\chat=-2$ can be derived by using \eqref{dwdc}: The fact that $\mu\searrow 4$ means $c\nearrow -2$ leads to
\begin{eqnarray*}
  \left|\frac{dw(c)}{dc}\right|    &=& O\left(\frac{1}{\sqrt{\mu-4}}\right) \qquad\mbox{(by Theorem \ref{Main_thm})}\\
            &=& O\left(\frac{2}{\sqrt{(\mu+2)(\mu-4)}}\right) \\
  &=& O\left(\frac{1}{\sqrt{-2-c}}\right) \qquad\mbox{(by the identity \eqref{cmu})}
\end{eqnarray*}
as $c\nearrow -2$, in which we have used $-z+1/2=w/\mu$ and the uniform boundedness of $J(q_c)$. (An estimate of the size of the Julia set is given in Lemma \ref{uniform_bound_q}.) 
\end{remark}

\begin{remark} \rm
The estimates in theorems \ref{Main_thm_2} and \ref{Main_thm} are optimal. 
For example, let $z=z(c)$ be the repelling fixed point of $q_c$ for $0\le c<1/4$, then 
$$
  \frac{dz(c)}{dc}=- \frac{1}{2\sqrt{1/4-c}}.
$$
(This means that the equality in \eqref{estimate2} is attained.)
Also,   let $z=z(\mu)$ be a pre-image point of  $1\in J(f_\mu)$, then 
$$
  \frac{dz(\mu)}{d\mu}=\pm \frac{1}{\mu\sqrt{\mu}\sqrt{\mu-4}}.
$$
\end{remark}

\section{Hausdorff convergence and dynamical degeneration}

In this section, we discuss the convergence of the Julia set $J(q_c)$ and the degeneration of dynamics of $q_c$ on $J(q_c)$ when $c$ approaches the boundary points $c=1/4$ or $-2$ of the Mandelbrot set. We also discuss corresponding properties for the map $f_\mu$ on its Julia set $J(f_\mu)$ by using \eqref{conjugacy} and \eqref{cmu}. The equality \eqref{cmu} shows that the map from real $\mu$ to real $c$ in $(-\infty, 1/4]$ is two-fold, therefore $\mu\searrow 1$ or $\mu\nearrow 1$ when $c\nearrow 1/4$, and $\mu\searrow 4$ or $\mu\nearrow -2$ as $c\nearrow -2$ (or more precisely, $\mu=1\pm\sqrt{4\epsilon}$ when $c=1/4-\epsilon$ and $\mu=1\pm\sqrt{9+4\epsilon}$ when $c=-2-\epsilon$, with $\epsilon>0$). It is easy to see that the dynamics of $f_{1-\sqrt{4\epsilon}}$ is a trivial copy of that of $f_{1+\sqrt{4\epsilon}}$, and the same triviality holds for $f_{1-\sqrt{9+4\epsilon}}$ and $f_{1+\sqrt{9+4\epsilon}}$. Hence, we shall restrict our discussion to the cases $\mu\searrow 1$ and $\mu\searrow 4$ only. 
 
\subsection{$c\nearrow 1/4$ or $\mu\searrow 1$}

Theorem \ref{Main_thm_2} imples one-sided H\"{o}lder continuity of the holomorphic motion as $c\nearrow 1/4$ along the real line:
The improper integral 
$$
z(1/4) := z(0)+ \lim_{\delta \searrow 0} \int_{0}^{1/4-\delta} \frac{dz(c)}{dc} ~dc
$$
exists in the Julia set $J(q_{1/4})$. 
In particular,  
\begin{equation}
|z(c)-z(1/4)| \le \sqrt{1/4-c} \label{Holderc14}
\end{equation}
for any real $c\in [0,1/4)$.

It is well-known that $q_c$ is hyperbolic for $c$ in the hyperbolic components of $\M$. Therefore, $(J(q_c), q_c)$, the restriction of $q_c$ to $J(q_c)$, is topologically conjugate to $(J(q_0), q_0)$ via a conjugacy $h_c(\cdot;0)$ from $J(q_0)$ to $J(q_c)$ for $c$ in the main cardioid of $\M$. Inequalities \eqref{estimate2} and \eqref{Holderc14} lead to a result that  the conjugacy $h_c(\cdot;0)$ 
converges uniformly to a semiconjugacy $h_{1/4}(\cdot;0):J(q_{0}) \to J(q_{1/4})$, $z(0)\mapsto z(1/4)$, as $c$ increases from $0$ to $1/4$. 
 As a matter of fact, $h_{1/4}(\cdot;0)$ is a conjugacy (see \cite{K} for example). So, $(J(q_{1/4}), q_{1/4})$ is topologically conjugate to $(\T, \mathcal{T})$, where 
$$\mathcal{T}:\T\to\T,~ t\mapsto 2t ~(\bmod~ 1)$$
is  the angle-doubling map.

\paragraph{Hausdorff distance.}
The estimate \eqref{Holderc14} and 
the existence of the conjugacy above between 
$(J(q_c),q_c)$ and $(J(q_{1/4}),q_{1/4})$ 
imply that the Hausdorff distance between 
$J(q_c)$ and $J(q_{1/4})$ is at most 
$\sqrt{1/4-c}$. 
The distance is also at least $\sqrt{1/4-c}$,
because the distance between the parabolic fixed point 
$1/2 \in J(q_{1/4})$ and the Julia set $J(q_c)$ 
is attained by the repelling fixed point $(1+\sqrt{1-4c})/2=1/2+\sqrt{1/4-c}$ of $q_c$. (Indeed, in the proof of Theorem \ref{Main_thm_2} we will show that $|z| \ge (1+\sqrt{1-4c})/2$ for any $z \in J(q_c)$.)
Hence we obtain the following:

\begin{cor}
For $0 \le c <1/4$, the Hausdorff distance between
the Julia sets $J(q_c)$ and $J(q_{1/4})$ is exactly $\sqrt{1/4-c}$.
\end{cor}

\begin{remark} \label{rk_unexpected} \rm
  From \eqref{dzdmu} and Theorem  \ref{Main_thm_2}, when $1<\mu<2$ we get 
\begin{eqnarray*}
  \left| \frac{dz(\mu)}{d\mu}\right| &\le& \frac{\mu-1}{2\mu}\cdot\frac{1}{\sqrt{1-4c}}+\frac{1+\sqrt{1+4c}}{2}\cdot\frac{1}{\mu^2} \qquad\mbox{(by Lemma \ref{uniform_bound_q})}\\
                                                  &=& \frac{1}{2\mu} +\frac{1+\sqrt{1+2\mu-\mu^2}}{2\mu^2} \qquad\mbox{(by the identity \eqref{cmu})} \\
                                                   &\le& \frac{2+\sqrt{2}}{2}.
\end{eqnarray*}
This unexpected result implies that the Hausdorff distance between $J(f_\mu)$ and $J(f_1)$ is at most $(2+\sqrt{2})(\mu-1)/2$ as $\mu\searrow 1$ along the real axis.
\end{remark}  

\paragraph{Notation.}
When variables $X, Y \ge 0$ satisfy 
$X/C \le Y \le CX$ with a uniform constant $C>1$, 
we denote this by $X \asymp Y$.

\begin{remark} \label{rkdim} \rm
  Let $\dim_H J(q_c)$ denote the Hausdorff dimension of $J(q_c)$. It has been known from \cite{HZ} that 
$$ 1<\dim_H J(q_{\frac{1}{4}})<\frac{3}{2}$$
 and that there exists $c_0<1/4$ such that for all $c\in [c_0,1/4)$ one has 
\begin{equation}
 \frac{d}{dc}\dim_H J(q_c)\asymp \left(\frac{1}{4}-c \right)^{\dim_H J(q_{\frac{1}{4}})-\frac{3}{2}}. \label{ddcdimH14}
\end{equation}
 Therefore, the derivative of the Hausdorff dimension of the Julia set $J(q_c)$ with respect to $c$ tends to infinity from the left of $1/4$, and the graph of $\dim_H J(q_c)$ versus $c$ has a vertical tangent on the left at $1/4$. The unexpected result in Remark \ref{rk_unexpected} suggests us to examine whether or not a similar vertical tangency holds for the graph of the Hausdorff dimension $\dim_H J(f_\mu)$ as $\mu$ decreases to  $1$. What we find is that the graph $\dim_HJ(f_\mu)$ versus $\mu$ has a horizontal tangent on the right at $1$. (See the proposition below.)
\end{remark}

\begin{prop}
 There exists $\mu_0>1$ such that 
$$ 
-\frac{d}{d\mu}\dim_H J(f_\mu)\asymp \left( \mu-1 \right)^{2\dim_H J(f_1)-2}
$$
for any $\mu\in (1,\mu_0]$.
\end{prop}
\paragraph{Proof.}
For $\mu\not=0$,  the affine map $G(\cdot, \mu)$ in \eqref{conjugacy} sends $J(f_\mu)$ to $J(q_{\mu(2-\mu)/4})$, thus
$$ \dim_H J(f_\mu)=\dim_H J(q_{\frac{\mu (2-\mu)}{4}})
$$
for $\mu\not=0$.
In particular, 
$$ \dim_H J(f_1)=\dim_H J(q_{\frac{1}{4}}).
$$
Consequently, for $\mu\not=0$ or $1$, 
\begin{eqnarray*}
   -\frac{d}{d\mu}\dim_H J(f_\mu)&=&\frac{d}{dc}\dim_H J(q_c)\cdot \frac{\mu-1}{2} \\
                                        &\asymp& \left( \frac{1}{4}-\frac{\mu(2-\mu)}{4}\right)^{\dim_HJ(f_1)-\frac{3}{2}}\cdot \frac{\mu-1}{2} \qquad\mbox{(by \eqref{ddcdimH14})}\\
                                                &=&\left(\frac{\mu-1}{2}\right)^{2\dim_H J(f_1)-2}\\
                                        &\asymp&\left(\mu-1\right)^{2\dim_H J(f_1)-2},
\end{eqnarray*}
as asserted. The value of $\mu_0$ can be obtained by solving $c_0=\mu_0(2-\mu_0)/4$ from \eqref{cmu} with the $c_0$ in Remark \ref{rkdim}.
\QED

\subsection{$c\nearrow -2$ or $\mu\searrow 4$}

For $c$ in the  exterior of $\M$, the restriction of $q_c$ to $J(q_c)$ is topologically conjugate to the one-sided left shift with two symbols. Therefore, there exists a conjugacy $h_c(\cdot;c_0)$ from $J(q_{c_0})$ to $J(q_c)$ for any fixed $c_0$ and $c$ not belonging to $\M$.
  The results \eqref{Main_estimate} and \eqref{eq2} give rise to a consequence that for 
 any semi-hyperbolic parameter $\chat \in \partial \M$, any parameter ray $\mathcal{R}(\theta)$ 
landing at $\chat$, and any $c_0\in\mathcal{R}(\theta)$, the conjugacy $h_c(\cdot;c_0)$ 
converges uniformly to a semiconjugacy 
\begin{equation}
  h_{\chat}(\cdot;c_0):J(q_{c_0}) \to J(q_\chat), \quad z(c)\mapsto z(\chat), \label{eq_semiconj}
\end{equation}
as $c \to \chat$ along $\mathcal{R}(\theta)$. 
This further implies that 
the Hausdorff distance between $J(q_c)$ and $J(q_\chat)$ 
is $O(\sqrt{|c-\chat|})$ as $c \to \chat$ along $\cR(\theta)$.

Let 
$$
\Sigma:=\Big\{
{\bf s}=\{s_0,s_1,s_2,\ldots\} |~ s_n=0\ \mbox{or}\ 1~ \mbox{for all}~ n\ge 0
\Big\}
$$ 
be the space consisting of sequences of $0$'s and $1$'s with the product topology, and  $\sigma$ be the left shift in $\Sigma$, $\sigma( {\bf s})={\bf s}'=\{s'_0, s'_1,s'_2,\cdots\}$ with $s'_i=s_{i+1}$.
 Fix $\theta\in\T-\{0\}$, the two points $\theta/2$ and $(\theta+1)/2$ divide $\T$ into two open semi-circles $\T_0^\theta$ and $\T_1^\theta$ with $\theta\in\T_0^\theta$. Let  $\theta$ be such an angle that $ \mathcal{T}^n(\theta)\not\in \left\{\frac{\theta}{2},\frac{\theta+1}{2}\right\}$ for all $n\ge 0$. Define the {\it kneading sequence} of $\theta$ under $\mathcal{T}$ as $\mathcal{E}(\theta)=\{ \mathcal{E}(\theta)_n \}_{n\ge 0}\in\Sigma$ with
\[  \mathcal{E}(\theta)_n=
\begin{cases}
 0 & \mbox{for}~ \mathcal{T}^n(\theta)\in \T_0^\theta \\
 1 & \mbox{for}~ \mathcal{T}^n(\theta)\in \T_1^\theta.
\end{cases}
\]
Note that the kneading sequence of non-recurrent $\theta$ is well-defined.

A point ${\bf e}\in\Sigma$ is said to be {\it aperiodic} if $\sigma^n({\bf e})\not={\bf e}$ for any $n\ge 0$. 
  Two points ${\bf a}$ and ${\bf s}$ in $\Sigma$ are said to be {\it equivalent} with respect to  aperiodic ${\bf e}\in \Sigma$, denoted by ${\bf a}\sim_{\bf e}{\bf s}$, if there is $k\ge 0$ such that $a_n=s_n$ for all $n\not= k$ and $\sigma^{k+1}({\bf a})=\sigma^{k+1}({\bf s})={\bf e}$.

In \cite{CK}, we proved that the semiconjugacy $h_\chat(\cdot; c_0)$ described in \eqref{eq_semiconj} leads to the following result. Let $\chat$ be a semi-hyperbolic parameter with an external angle $\theta$ and ${\bf e}=\mathcal{E}(\theta)$ be the kneading sequence of $\theta$. Then $(J(q_\chat), q_\chat)$ is topologically conjugate to $(\Sigma/{\sim_{\bf e}}, \tilde\sigma)$, where $\tilde\sigma$ is induced by the shift transformation $\sigma$.

The dynamical degeneration  for $q_c$ as $c\nearrow -2$ along  real axis, namely, along $\mathcal{R}(\theta)$ with $\theta=1/2$  is the same as the one for the logistic map 
$f_\mu$ as $\mu\searrow 4$ along the  real axis. Since  $\mathcal{E}(1/2)=\{0,1,1,\ldots\}$, the dynamical degeneration of $q_c$ at $c=-2$ or $f_\mu$ at $\mu=4$ is that both $(J(q_{-2}),q_{-2})$ and $(J(f_4),f_4)$ are topologically  conjugate to $(\Sigma/{\sim_{\{0,1,1,\ldots\}}},\tilde\sigma)$.

Note that there is another way to interpret the kneading sequence:
 For the  logistic map $f_\mu$, the Julia set $J(f_\mu)$  is a Cantor set contained in the real interval $[0,1]$ when $\mu$ is real and greater than $4$. When $\mu = 4$, $J(f_4)$ is the whole interval $[0,1]$.
 If the critical point $1/2$ belongs to the Julia set, one can define the  {\it kneading sequence}  $I(f_\mu)=\{I(f_\mu)_n\}_{n\ge 0}$ for $f_\mu$ by $I(f_\mu)_n=1$ if $f_\mu^{1+n}(1/2)\in [0,1/2] \cap J(f_\mu)$ and  $I(f_\mu)_n=0$ if $f_\mu^{1+n}(1/2)\in [1/2, 1] \cap J(f_\mu)$. The sequence $I(f_\mu)$ is well-defined if $f_\mu^{1+n}(1/2)\not=1/2$ for all $n\ge 0$.  Then, it is not difficult to see that $\mathcal{E}(1/2)=I(f_4)$.

\section{Proof of Theorem \ref{Main_thm_2}}

Assume $0\le c < 1/4$. Let $r=(1+\sqrt{1-4c})/2$, the distance of the repelling fixed point of $q_c$ from the critical point $0$. 
Suppose $z\in J(q_c)$ and $|z|<r$. Then, 
$|q_c(z)|=|z^2+c|<r^2+c=r$. Thus, $|q_c^n(z)|<r$ for all $n\ge 0$. Since this is an open condition, 
this means that $z$ belongs to the filled Julia set but not to the Julia set, a contradiction. 
Therefore, 
 $$ 
  \inf_{z\in J(q_c)}|z|=r,
$$
and 
$|Dq_c(z)|=2|z|\ge 1+\sqrt{1-4c}$ for all $z\in J(q_c)$ (where $D=d/dz$). For any $c$ not belonging to $\M$ or in a hyperbolic component of $\M$, and for any $z =z(c) \in J(q_c)$, in \cite{CK} we proved the following derivative formula
$$
\frac{dz(c)}{dc} = 
-\sum_{n = 1}^\infty\frac{1}{Dq_c^{{n}}({z}(c))}.
$$
Hence, 
$$
  \left|\frac{dz(c)}{dc}\right| \le \sum_{n\ge 1} \frac{1}{|Dq_c^n(z)|} \le \sum_{n\ge 1}\frac{1}{\left(1+\sqrt{1-4c}\right)^n}=\frac{1}{2\sqrt{1/4-c}}.
$$
\QED

\begin{remark}\rm ~ \\
(i) Compared with our proof of Theorem \ref{thm_Main_estimate} in \cite{CK} for $c$ approaching a semi-hyperbolic parameter, or even with that of Theorem \ref{Main_thm}, to come in the next section, we find that the proof of Theorem \ref{Main_thm_2} is surprisingly simple. \\
(ii) In fact, by combining with the lemma below, the  Julia set $J(q_c)$  locates inside the annulus $\{z|~ (1+\sqrt{1-4c})/2\le|z|\le (1+\sqrt{1+4c})/2\}$ for $0\le c<1/4$. 
\end{remark}

\begin{lem}  \label{uniform_bound_q} 
 $|z|\le (1+\sqrt{1+4|c|})/2$ for any  $z\in J(q_c)$ and  $c\in\mathbb{C}$. 
\end{lem}
\paragraph{Proof.}
 Let $M=M(c)= (1+\sqrt{1+4|c|})/2$ and notice that $M^2-M-|c|=0$. If $z$ is such a point that $|z|=M+R$ for some $R>0$, then $|q_c(z)|-|z|=|z^2+c|-|z|\ge |z|^2-|z|-|c|=(2M-1)R+R^2 $. This implies that the orbit of $z$ tends to infinity thus $z\not\in J(q_c)$.
\QED

\section{Proof of Theorem \ref{Main_thm}}

Note that $f_\mu(1) = 0$ and the origin $z = 0$ is a repelling fixed point of multiplier $\mu$ for any real $\mu \ge 4$. Hence there exists a linearizing coordinate
$\phi_\mu:\tilde{U}_0 \to \C$ defined on 
a fixed  neighborhood $\tilde{U}_0$ of $0$ such that 
\begin{enumerate} 
\item
for any $z \in f_\mu^{-1}(\tilde{U}_0)\cap \tilde{U}_0$,
$$
\phi_\mu(f_\mu(z)) = \mu \phi_\mu(z).
$$
\item
$ |\phi_\mu(z)| \asymp |z|$ when $\mu$ is sufficiently close to $4$.
\end{enumerate}

Fix a point $z = z(\mu)$ in the Julia set $J(f_\mu)$ for $\mu >4$.
Set $f := f_\mu$ and $z_n =z_n(\mu): = f_\mu^n(z)~(n \ge 0)$. 
As in \cite{CKLY}, we can derive 
$$
\frac{dz_{n + 1}}{d\mu}
 = 
Df(z_n)\frac{dz_{n}}{d\mu}
 + \frac{1}{\mu} z_{n + 1} 
$$
(where $D = d/dz$) and thus we have a `formal' expansion
\begin{equation}
\frac{dz}{d\mu}
 =  
 -\frac{1}{\mu} \sum_{n \ge 1}\frac{z_{n}}{Df^n(z)}. \label{formal_derivative}
\end{equation}
Now suppose that $\mu$ is sufficiently close to $4$. 
We will show that the formal expansion above converges absolutely.

The main idea is to consider a fixed `singular' metric of the following explicit form
$$
\gamma(z) |dz|:=  \frac{|dz|}{\sqrt{|z|~|z-1|}}
$$
on $\C-\brac{0,1}$ (inspired by \cite{Ka}; see also \cite{Ro}). 
Let us fix a $z \in J(f)\cap (0,1)$. 
Then it is easy to see $\gamma(z) \ge 2$. 
Moreover, 
$$
\frac{\gamma(f(z))|Df(z)|}
{\gamma(z)} = 
\frac{2 \sqrt{\mu}~|z-1/2|}
{\sqrt{1-f(z)}}.
$$
Note that $f(z) \in J(f) \subset [0,1]$, hence $1-f(z) \ge 0$. By $f(z) = \mu/4-\mu(z-1/2)^2$,
 we have 
$$
\frac{\gamma(f(z))|Df(z)|}
{\gamma(z)} = 
\frac{2 \sqrt{\mu/4 -f(z)}}
{\sqrt{1-f(z)}}.
$$
The right hand side takes its infimum as $f(z) \searrow 0$:
$$
\frac{\gamma(f(z))|Df(z)|}
{\gamma(z)} \ge 
2 \sqrt{\frac{\mu}{4}}
 = \sqrt{\mu}
  =:A ~(= 2 + O(\mu-4)).
$$
Hence for any fixed $z= z_0 \in J(f)$
whose forward orbit never lands on the fixed point $0$, with the help of the  identity
$$
\frac{\gamma(z_n)}{\gamma(z)}Df^n(z)=\prod_{k=0}^{n-1} \frac{\gamma(z_{k+1})}{\gamma(z_k)}Df(z_k),
$$
we obtain
\begin{equation}
\frac{1}{|Df^n(z)|} \le \frac{\gamma(z_n)}{A^n \gamma(z)}
\le \frac{\gamma(z_n)}{2A^n}
\asymp \frac{\gamma(z_n)}{A^n}. \label{eq0}
\end{equation}
Now we have 
\begin{equation}\label{eq01}
\abs{\frac{dz}{d\mu}}
 \le 
 \mu^{-1} \sum_{n \ge 1}\frac{|z_{n}|}{|Df^n(z)|}
 \le 
 \mu^{-1} \sum_{n \ge 1}\frac{|z_{n}|\gamma(z_n)}{2A^n}
 \asymp
\sum_{n \ge 1}\frac{|z_{n}|\gamma(z_n)}{A^n}
\end{equation}
for such a $z$.
This implies that 
if the forward orbit $\brac{z_n}_{n \ge 0}$ is 
a certain distance away from $0$ (and $1$), 
then $|z_n|\gamma(z_n) \asymp 1$ and 
the derivative $\dfrac{dz}{d\mu}$ is uniformly bounded
by a constant. More precisely, we have:
\begin{prop}
For any real $\mu \ge 4$ and any forward orbit $\brac{z_n}_{n \ge 0}$, 
if there exists some $\delta >0$ such that 
$
\delta \le z_n  \le 1-\delta
$
for all $n$, then 
$$
\abs{\dfrac{dz}{d\mu}} 
\le \frac{1}{8\delta}.
$$
\end{prop}

\paragraph{Proof.}
By assumption we have $|z_n| \le 1$ and $\gamma(z_n) \le 1/\delta$.
Thus (\ref{eq01}) implies
\begin{align*}
\abs{\frac{dz}{d\mu}}
 \le 
 \mu^{-1} \sum_{n \ge 1}\frac{|z_{n}|\gamma(z_n)}{2A^n}
\le
\mu^{-1}\cdot \frac{1}{2\delta}\cdot\frac{A^{-1}}{1-A^{-1}}.
\end{align*}
Since $\mu = A^2 \ge 4$ we have the desired estimate. 
\QED

\paragraph{Non-pre-fixed case.}
Next we suppose that $\brac{z_n}_{n \ge 0}$ accumulates on $0$ but never lands on $0$. 
We may assume that $f^{-1}(\tilde{U}_0)$ is the union of disjoint 
neighborhoods $U_0$ of $0$ and $U_1$ of $1$. 

Now there exist $N$ and $m \ge 1$ such that 
\begin{itemize}
\item
$z_{N-1} \in U_1 \cap J(f)$,
\item 
$z_{N},z_{N + 1},\, \ldots, \, z_{N + m-1} \in U_0 \cap J(f)$,
\item 
$z_{N + m} \in (\tilde{U}_0-U_0) \cap J(f)$.
\end{itemize}
By the linearizing coordinate $\phi = \phi_\mu:\tilde{U}_0 \to \C$, 
we have
$$
\phi(z_{N + m}) = \phi(f^{m-i}(z_{N + i})) 
= \mu^{m-i}\phi(z_{N + i})
$$
for $0 \le i \le m$. 
Since $|\phi(z)| \asymp |z|$,  we have 
$|\phi(z_{N + m })| \asymp |z_{N + m }| \asymp 1$ and thus
$$
|z_{N + i}| \asymp \mu^{-m + i}.
$$
Since $\gamma(z) \asymp 1/\sqrt{z}$ when $0<z<1/2$, we have
$$
\sum_{i = 0}^{m}
\frac{|z_{N + i}|\gamma(z_{N + i})}{A^{N + i}}
\asymp 
\sum_{i = 0}^{m}
\frac{\mu^{-m + i}\cdot \sqrt{\mu^{m - i}}}{A^{N + i}}
 = \frac{m+1}{A^{N + m}}.
$$
Hence the sum for `near zero' orbit points 
$z_{N}, \ldots, z_{N + m} \in \tilde{U}_0$ 
are bounded by $(m + 1)/A^{N + m}$.
Since $(m + 1)/A^m \to 0$ as $m \to \infty$, 
there exits a constant $K>0$ independent of $m$ with 
$(m + 1)/A^{N + m} \le  K/A^N$.
\footnote{By calculating the function $z \mapsto (z + 1)/A^z$,
 we will find $K = A/(e\log A) = 2/(e\log 2) + O(\mu-4)$.}

Let us give an estimate of $\dfrac{|z_{N -1}|}{|Df^{N-1}(z)|}$
for $z_{N-1} \in U_1$. 
When $N = 1$ this term is not counted in the formal sum expansion of $dz/d\mu$, so we may assume that $N \ge 2$.

One can show that $1-z_{N-1} \asymp \mu^{-m-1}$. 
Since 
$$
|Df(z)| = \mu ~ |1-2z| = 2\sqrt{\mu} \sqrt{\frac{\mu}{4} -f(z)},
$$
we have 
$$
|Df(z_{N-2})|
= 2\sqrt{\mu} \sqrt{\frac{\mu-4}{4} + 1-z_{N-1}} 
\ge \sqrt{\mu} \sqrt{\mu-4}.
$$
Moreover, by using \eqref{eq0} and $z_{N-2} \approx 1/2$, we have 
$$
\frac{1}{|Df^{N-2}(z)|} 
 = O\paren{\frac{\gamma(z_{N-2})}{A^{N-2}}} 
 = O\paren{\frac{1}{A^N}}.
$$
Hence
$$
\frac{|z_{N -1}|}{|Df^{N-1}(z)|}
 = 
\frac{|z_{N -1}|}{|Df^{N-2}(z)\cdot Df(z_{N-2})|}
= O\paren{\frac{1}{A^{N}} 
\cdot \frac{1}{\sqrt{\mu-4}}}.
$$

By assumption, we have infinitely many `near singular' orbit points
$$
z_{N_1-1}, \ldots, z_{N_1 + m_1}, 
z_{N_2-1}, \ldots, z_{N_2 + m_2}, 
\ldots \in \tilde{U}_0 \cup U_1
$$
with strictly increasing $N_j$.
Hence the original sum \eqref{formal_derivative} can be estimated as follows:
\begin{align}
\abs{\frac{dz}{d\mu}}
&\le \sum_{z_n \notin \tilde{U}_0 \cup U_1}\frac{|z_{n}|\gamma(z_n)}{A^n}
 + 
\sum_{z_n \in \tilde{U}_0 }\frac{|z_{n}|\gamma(z_n)}{A^n}
 + \sum_{z_n \in U_1}\frac{|z_{n}|\gamma(z_n)}{A^n}
\label{eq1}\\
& \asymp
\sum_{z_n \notin \tilde{U}_0 \cup U_1}\frac{1}{A^n}
 + 
\sum_{j\ge 1}
\frac{m_j + 1}{A^{N_j + m_j}}
+ 
\sum_{j\ge 1}
O\paren{\frac{1}{A^{N_j}} 
\cdot \frac{1}{\sqrt{\mu-4}}}
\nonumber\\
& 
= 
O(1)
 + 
\sum_{j\ge 1}\frac{K}{A^{N_j}}
 + 
O\paren{\frac{1}{\sqrt{\mu-4}}}
\cdot \sum_{j\ge 1}\frac{1}{A^{N_j}}\nonumber\\
& = O\paren{\frac{1}{\sqrt{\mu-4}}}
~~(\mu \searrow 4). \nonumber
\end{align}
Since the diameter of $\tilde{U}_0$ is fixed for $\mu \approx 4$,
 one can check that the estimate above is independent of both 
 $\mu$ and $z = z(\mu) \in J(f_\mu)$.

\paragraph{Pre-fixed case.}
Next we consider the case where the orbit of $z = z(\mu)$
eventually lands on $0$. 
Assume that $z_{n} = 0$ if and only if $n \ge N \ge 0$.
Then the derivative is written as a finite sum
\begin{align*}
\frac{dz}{d\mu} &= -\frac{1}{\mu} \sum_{1 \le n < N} \frac{z_n}{Df^n(z)}
\end{align*}
and the sum is divided into 
$$
\sum_{1 \le n < N} = 
\sum_{z_n \notin \tilde{U}_0 \cup U_1}
 + 
\sum_{z_n \in \tilde{U}_0 } 
 + \sum_{z_n \in U_1}
$$
as \eqref{eq1} in the previous case. 
Hence we obtain the same estimate $|dz/d\mu| =  O(1/\sqrt{\mu-4})$
without extra effort. 
\QED

\section*{Acknowledgments}
Chen was partly supported by MOST 106-2115-M-001-007.
Kawahira was partly supported by JSPS KAKENHI Grant Number 16K05193. 
They thank the hospitality of Academia Sinica, Nagoya University, 
RIMS in Kyoto University, and Tokyo Institute of Technology
where parts of this research were carried out.

\vspace{1cm}

{~}\\
Yi-Chiuan Chen \\
Institute of Mathematics\\
Academia Sinica\\
Taipei 10617, Taiwan\\
YCChen@math.sinica.edu.tw

\vspace{.5cm}

{~}\\
Tomoki Kawahira\\
Department of Mathematics\\ 
Tokyo Institute of Technology\\
Tokyo 152-8551, Japan\\
kawahira@math.titech.ac.jp  \\

\noindent
Mathematical Science Team\\ 
RIKEN Center for Advanced Intelligence Project (AIP)\\
1-4-1 Nihonbashi, Chuo-ku\\ 
Tokyo 103-0027, Japan

\end{document}